\documentstyle[12pt]{amsart}
\newtheorem{theorem}{Theorem}[section]
\newtheorem{corollary}[theorem]{Corollary}

\newtheorem{lemma}[theorem]{Lemma}

\numberwithin{equation}{section}

\begin{document}

\title [An $\varepsilon$-regularity Theorem For The Mean Curvature Flow]
 {An $\varepsilon$-regularity Theorem For The Mean Curvature Flow}
\author{Xiaoli Han, Jun Sun}

\address{Department of Mathematical sciences, Tsinghua university, Beijing, 100084,
P. R. China.} \email{xlhan@@math.tsinghua.edu.cn}

\address{Math. Group, The Abdus Salam International Centre for Theoretical Physics\\ Trieste 34100, Italy.}
\email{jsun@@ictp.it}


\date{}

\maketitle

\begin{abstract}
In this paper, we will derive a small energy regularity theorem for the mean curvature flow of arbitrary dimension and codimension. It says that
if the parabolic integral of $|A|^2$ around a point in space-time is small, then the mean curvature flow cannot develop singularity
at this point. As an application, we can prove that the 2-dimensional Hausdorff measure of the singular set of the mean curvature flow
from a surface to a Riemannian manifold must be zero.
\end{abstract}


\section{Introduction}

\allowdisplaybreaks

Recently years, the geometric flow becomes a study focus in geometric analysis. Among several basic problems, the regularity of the geometric flow
attracts many attentions of the mathematicians. Roughly speaking, there are two problems that people concern.

On one hand, they want to know whether,
or more precisely, under what condition, the flow can exist globally and converge to an ideal object. For example, in \cite{ES}, Eells-Sampson
proved that if the sectional curvature of the target manifold is nonpositive, then the harmonic map heat flow would converge to a harmonic map at infinity.
For mean curvature flow, Ecker-Huisken (\cite{EH}), Chen-Li-Tian (\cite{CLT}), Wang (\cite{Wa2}) proved some global existence and convergence results
for the graphic case.

On the other hand, as singularity always occurs, one wants to study the properties of the singularity, such as, what does the possible singularity look
like, and what is the behavior of singular set, etc. An important tool to study this is monotonicity formula. Using the monotonicity formula for mean curvature flow, Huisken (\cite{Hu1}) studied the Type I singularity of mean convex mean curvature flow. Chen-Li (\cite{CL1}) and Wang (\cite{Wa1}) independently
showed that there is no Type I singularity for symplectic mean curvature flow. Another important tool is small energy regularity theorem. It describes the local
behavior around the singular point.

\vspace{.2in}

Using a blow-up analysis and mean value inequality, K. Ecker (\cite{Ec}) proved the following small energy regularity
theorem:

\begin{theorem}\label{thm1.1}(\cite{Ec})
Suppose a family of surfaces $(M_t)_{t\in (0,T]}$ moves by their mean curvature in $\textbf{R}^3$. Then there exist
constants $\varepsilon_0>0$ and $c_0>0$ such that for $x_0\in \textbf{R}^3$ and $\rho\in (0,\sqrt{T}]$, the inequality
\begin{equation*}
\sup_{[T-\rho^2,T]}\int_{M_t\cap B_{\rho}(x_0)}|A|^2 \leq \varepsilon_0
\end{equation*}
implies the estimate
\begin{equation*}
\max_{\sigma\in [0,\rho]}\sigma^2\sup_{t\in[T-(\rho-\sigma)^2,T)}\sup_{M_t\cap B_{\rho-\sigma}(x_0)}|A|^2 \leq c_0.
\end{equation*}
\end{theorem}

\vspace{.2in}

Later on, T. Ilmanen (\cite{Il}) generalized Ecker's result to arbitrary dimension and codimension and the assumption is weaken
to a parabolic one:

\begin{theorem}\label{thm1.2}(\cite{Il})
There is a constant $\varepsilon_0=\varepsilon_0(n,k)$ such that if $(M_t^k)_{t\in (0,1]}$ is a mean curvature flow smoothly
immersed in $B_1\subset \textbf{R}^n$, and
\begin{equation*}
r^{-k}\int^{t}_{t-r^2}\int_{M_t\cap B_{r}(x)}|A|^2d\mu_t dt\leq \varepsilon^2 \leq \varepsilon_0^2
\end{equation*}
whenever $B_r(x)\times [t-r^2,t) \subset B_1 \times [0,1)$, then $M_t \cap B_1$ can be smoothly extended to $B_1 \times \{1\}$ and
\begin{equation*}
|A(x,t)|^2 \leq C(n,k)\max \left(\frac{1}{1-|x|},\frac{1}{t^{\frac{1}{2}}}\right)\varepsilon
\end{equation*}
for all $(x,t)\in B_1 \times (0,1]$.
\end{theorem}

In this paper, we will derive another $\varepsilon$-regularity theorem for the mean curvature flow for arbitrary dimension and codimension. Unlike
the proof of the above theorems, we obtain the result by combining a blow up argument and a $\kappa$-noncollapsing theorem. Our main theorem is as follows:

\vspace{.2in}

\noindent {\bf Main Theorem } {\it Let $M$ be an $\bar{n}$-dimensional complete Riemannian manifold
and $\Sigma_0$ be an n-dimensional submanifold in $M$. Suppose $\Sigma_0$ evolves by the mean curvature flow in $M$.
Then there exist constants $\varepsilon_0, r$
depending on $\Sigma_0$ and $M$, such that if for all $0\leq \rho <\rho' \leq \frac{r}{2}$
\begin{equation*}
((\rho')^2-\rho^2)^{-\frac{n}{2}}\int^{T-\rho^2}_{T-(\rho')^2}\int_{\Sigma_t \cap B_{((\rho')^2-\rho^2)^{\frac{1}{2}}(X_0)}}|A|^2d\mu_t dt \leq \varepsilon_0,
\end{equation*}
then we have
\begin{equation*}
\max_{\sigma\in (0,\frac{r}{2}]}\sigma^2\sup_{t\in[T-(r-\sigma)^2,T-(\frac{r}{2})^2]}\sup_{\Sigma_t\cap B_{r-\sigma}(X_0)}|A|^2 \leq C,
\end{equation*}
where $C$ depends on $\varepsilon_0$, $\Sigma_0$ and $M$.
}

\vspace{.2in}

Our Main Theorem says that if the parabolic integral of $|A|^2$ around a point $(X_0,T)$ is small, then $(X_0,T)$ cannot be a singular point.
Note that our assumption is weaker than both of Theorem \ref{thm1.1} and Theorem\ref{thm1.2}.

\vspace{.2in}

An application of $\varepsilon$-regularity theorem is to estimate the size of singular set. For example, Schoen-Uhlenbeck (\cite{SchU})
used the $\varepsilon$-regularity theorem for energy minimizer map to prove that the $m-2$-dimensional Hausdorff measure of the singular set of
an E-minimizer is 0. Using Theorem \ref{thm1.1}, Ecker proved similar result for mean curvature flow.

Applying our Main Theorem to the case $n=2$, we can get the following result:

\begin{corollary}\label{cor1.3}
 Let $\Sigma_0$ be a surface moving along its mean curvature in an $\bar{n}$-dimensional complete Riemannian manifold $M$ and $\mathcal F$ be the singular
set. Then ${\cal H}^2 ({\mathcal F})=0.$
\end{corollary}

\vspace{.1in}

\emph{Acknowledgement: The first author was supported by NSF in China (grant No.: 1090 1088). The second author thanks Abdus Salam ICTP for the good work condition. The authors would like to thank Professor Jiayu Li for his helpful discussion.}

\section{Proofs of The Main Theorem}

\vspace{.1in}

\textbf{Proof of The Main Theorem:}
We prove it by contradiction. Suppose the conclusion is false, then
$$\limsup_{r\to 0}\max_{\sigma\in (0,\frac{r}{2}]}\sigma^2\sup_{t\in[T-(r-\sigma)^2,T-(\frac{r}{2})^2]}\sup_{\Sigma_t\cap B_{r-\sigma}(X_0)}|A|^2
=\infty.$$
Then there exists a sequence $\{r_k\}$ with $r_k\to 0$ such that
$$\max_{\sigma\in (0, r_k/2]}\sigma^2 \max_{[T-(r_k-\sigma)^2,
T-(r_k/2)^2]} \max_{\Sigma_t\cap B_{r_k-\sigma}(X_0)}|A|^2\\
\rightarrow +\infty.
$$
We choose $\sigma_k\in (0, r_k/2]$ such that
$$\sigma_k^2 \max_{[T-(r_k-\sigma_k)^2, T-(r_k/2)^2]}
\max_{\Sigma_t\cap B_{r_k-\sigma_k}(X_0)}|A|^2=\max_{\sigma\in (0,
r_k/2]}\sigma^2 \max_{[T-(r_k-\sigma)^2, T-(r_k/2)^2]}
\max_{\Sigma_t\cap B_{r_k-\sigma}(X_0)}|A|^2.
$$
Let $t_k\in [T-(r_k-\sigma_k)^2, T-(r_k/2)^2]$ and
$F(x_k,t_k)=X_k\in\bar{B}_{r_k-\sigma_k}(X_0)$ satisfying
$$\lambda_k^2=|A|^2(X_k)=|A|^2(x_k, t_k)=\max_{[T-(r_k-\sigma_k)^2, T-(r_k/2)^2]}
\max_{\Sigma_t\cap B_{r_k-\sigma_k}(X_0)}|A|^2.$$ Obviously, we have
$(X_k, t_k)\to (X_0, T)$ and $\lambda_k^2\sigma_k^2\to\infty$. In
particular,
\begin{equation}\label{e0.2} \max_{[T-(r_k-\sigma_k/2)^2, T-(r_k/2)^2]} \max_{\Sigma_t\cap
B_{r_k-\sigma_k/2}(X_0)}|A|^2\leq 4\lambda_k^2, \end{equation} and
hence \begin{equation}\label{e0.3}\max_{[t_k-(\sigma_k/2)^2, t_k]}
\max_{\Sigma_t\cap B_{r_k-\sigma_k/2}(X_0)}|A|^2\leq 4\lambda_k^2.
\end{equation}

We choose a normal coordinates in $B_r(X_0)$ using the exponential
map, where $B_r(X_0)$ is a metric ball in $M$ centered at $X_0$ with
radius $r$ ($0<r<i_M/2$). We express $F$ in its coordinates
functions. Consider the following sequences,
\begin{equation}\label{e0.6} F_k(x, s)=\lambda_k(F(x_k+x,
t_k+\lambda_k^{-2}s)-F(x_k, t_k)),~~~~~~ s\in
[-\lambda_k^2\sigma_k^2/4, \lambda_k^2(T-t_k)].\end{equation} We
denote the rescaled surfaces by $\Sigma^k_s$ in which $d\mu^k_s$ is
the induced area element from $M$.  Therefore,
\begin{eqnarray*}
|A_k|^2 &=&\lambda_k^{-2}|A|^2,\\
H_k&=&\lambda_k^{-1}H,\\ |H_k|^2 &=&\lambda_k^{-2}|H|^2.
\end{eqnarray*}

Set $t=t_k+\lambda_k^{-2}s$, it is easy to check that
\begin{eqnarray*}
\frac{\partial F_k}{\partial s}&=&\lambda_k^{-1}\frac{\partial
F}{\partial t}.
\end{eqnarray*}

Therefore, it follows that the rescaled surface also evolves by a mean
curvature flow
\begin{equation}\label{e0.1} \frac{\partial F_k}{\partial s}=H_k
\end{equation} in $B_{\lambda_k\sigma_k}(0)$, where $s\in [-\lambda_k^2\sigma_k^2/4,
\lambda_k^2(T-t_k)]. $

By (\ref{e0.2}) and (\ref{e0.3}) we see that,
$$|A_k|(0, 0)=1,~~~~~~~~~|A_k|^2\leq 4$$ in $B_{\lambda_k\sigma_k}(0)$
and $s\in [-\lambda_k^2\sigma_k^2/4, 0]. $ Since we have
$\lambda_k^2\sigma_k^2\to\infty$,
thus by Arzela-Ascoli theorem, $\Sigma^k_s\to\Sigma^\infty_s$ in
$C^2(B_R(0)\times [-R, 0])$ for any $R>0$ and any $B_R(0)\subset
{\mathbb {R}}^{\bar{n}}$. By (\ref{e0.6}), we know that $\Sigma^\infty_s$ is
defined on $(-\infty, 0]$.  $\Sigma^\infty_s$ also evolves
along the mean curvature flow in ${\mathbb{R}}^{\bar{n}}$ with the Euclidean
metric and
$$|A_\infty|(0, 0)=1,~~~~~~~~~|A_\infty|^2\leq 4.$$

\vspace{.2in}

In order to proceeding further, we need some preparations. Let us first recall a $\kappa$-noncollapsing theorem proved by Chen-Yin (\cite{CY}):

\begin{theorem}\label{thm3.1}(\cite{CY})
Let $(M^{\bar{n}}, \bar{g})$ be a complete Riemannian manifold of dimension
$\bar{n}$ with bounded curvature and the injectivity radius is bounded from below
by a positive constant, i.e., there are constants $\bar{C}$ and $\bar{\delta}$ such that
\begin{equation}\label{e3.1}
|\bar{R}m|(x)\leq \bar{C}, \ \ and  \ \ inj(M^{\bar{n}}, x)\geq \bar{\delta} >0, \ \ for \ all \ x\in M^{\bar{n}}.
\end{equation}
Let $X: \Sigma^n \to M^{\bar{n}}$ be a complete isometrically immersed manifold
with bounded second fundamental form $|h^{\alpha}_{ij}|\leq C$ in $M^{\bar{n}}$,
then there is a positive constant $\delta=\delta(\bar{C},\bar{\delta},C,\bar{n})$
such that the injectivity radius of $\Sigma^n$ satisfies
\begin{equation}\label{e3.2}
inj(\Sigma^n, x)\geq \delta >0, \ \ for \ all \ x\in\Sigma^n.
\end{equation}
\end{theorem}

As the second fundamental form and its derivatives are all uniformly bounded
on $\Sigma^{\infty}_{s}$, by the same argument as in the proof of Lemma 4.1 in \cite{HanS2},
we can obtain the following $\kappa$-noncollapsing theorem. Note that for our case, the ambient space
is just the Euclidean space $\textbf{R}^{\bar{n}}$.

\begin{lemma}\label{lem3.2}
There exist constants $\kappa_{0}=\kappa_{0}(n)$ and $r_{0}=r_{0}(n,\bar{n})$ such that
$Vol_s(B_s(p,\rho))\geq \kappa_0 \rho^n$ as long as $\rho\leq r_{0}$
and $B_s(p, \rho)\subset \Sigma^{\infty}_s$, $s\in (-\infty,0]$.
\end{lemma}

\textbf{Proof:} Using Theorem \ref{thm3.1} we know that there
exists a constant $\iota=\iota(\bar{n})$
such that
\begin{equation}\label{e5}
inj(\Sigma_{s}^{\infty})\geq \iota, \ \ for \ s\in (-\infty,0].
\end{equation}
Moreover, by Gauss equation we have,
\begin{equation*}
\max_{\Sigma_{s}^{\infty}}|Rm|\leq C_1=C(n).
\end{equation*}
Then by volume comparison theorem, for $\rho \leq\iota$, $p\in
\Sigma_{s}^{\infty}$, we have
\begin{equation*}
Vol_s(B_{s}(p,\rho))\geq \bar{V}(\rho,C_{1},n),
\end{equation*}
where $\bar{V}(\rho,C_{1},n)$ is the volume of the geodesic ball of
radius $\rho$ in n-dimensional space form with constant curvature
$C_{1}$. But it is well known that
\begin{eqnarray*}
\bar{V}(\rho,C_{1},n)
         & =&   \int_{0}^{\rho}\int_{S^{n-1}}(\frac{\sin(\sqrt{C_{1}}r)}{\sqrt{C_{1}}})^{n-1}dS_{n-1}dr\nonumber\\
         & =&   n\omega_{n}\int_{0}^{\rho}(\frac{\sin(\sqrt{C_{1}}r)}{\sqrt{C_{1}}})^{n-1}dr\nonumber\\
         & =&   \frac{n\omega_{n}}{C_{1}^{\frac{n}{2}}}\int_{0}^{\sqrt{C_{1}}\rho}\sin^{n-1} r dr.
\end{eqnarray*}
On the other hand, $\lim_{r\to 0}\frac{\sin r}{r}=1$. Thus we can
choose
$r_{0}=r_{0}(C_{1},\iota)=r_{0}(n,\bar{n})<\iota$
such that for all $r\leq \sqrt{C_{1}}\rho \leq \sqrt{C_{1}}r_{0}$,
we have
$$\frac{\sin r}{r} \geq \frac{1}{2}.$$
Then for $\rho\leq r_{0}$,
\begin{equation*}
\bar{V}(\rho,C_{1},n)
   \geq \frac{n\omega_{n}}{C_{1}^{\frac{n}{2}}}\int_{0}^{\sqrt{C_{1}}\rho}\frac{r^{n-1}}{2^{n-1}} dr
   = \frac{\omega_{n}}{2^{n-1}}\rho^{n},
\end{equation*}
i.e.,
\begin{equation*}
Vol_s(B_{s}(p,\rho))\geq \frac{\omega_{n}}{2^{n-1}}\rho^{n},  \ \ for
\ \rho\leq r_{0}, \ s\in (-\infty,0].
\end{equation*}
This proves Lemma \ref{lem3.2}. \hfill Q.E.D.

\vspace{.2in}

As for the blow-up flow, the derivatives of the second fundamental form are all uniformly bounded,
we know from the evolution equations that there exist positive constants $C_2$, $C_3$ and $C_4$, such that
\begin{equation}\label{e3.3}
|\nabla A| \leq C_2 \ \ \ on \ \ \  \Sigma^{\infty}_{s}, \ \ for  \ s\in (-\infty,0],
\end{equation}
\begin{equation}\label{e3.4}
\left|\frac{\partial}{\partial s}|A|\right| \leq C_3 \ \ on \ \ \ \Sigma^{\infty}_{s}, \ \ for  \ s\in (-\infty,0],
\end{equation}
and
\begin{equation}\label{e3.5}
\left|\frac{\partial}{\partial s}g\right| \leq C_4 \ \ \ on \ \ \ \Sigma^{\infty}_{s}, \ \ for  \ s\in (-\infty,0].
\end{equation}

Using the $\kappa$-noncollapsing theorem, we can obtain the following
$\varepsilon_0$-regularity theorem:

\begin{theorem}\label{thm3.3}
For any $0<\varepsilon\leq \min\{r_0^{n+3},\frac{1}{2^{n+3}}\}$  if
$\int^{0}_{-1}\int_{\Sigma^{\infty}_{s}\cap B_{1}(0)}|A|^2d\mu_s^{\infty}ds\leq \varepsilon$, then
we have
\begin{equation}\label{e3.6}
\max_{-1\leq s \leq 0}\max_{\Sigma^{\infty}_{s}\cap B_{\frac{1}{2}}(0)}|A| \leq
\left(\sqrt{\frac{1}{\kappa_0}}+C_2+C_3+C_2C_4\right)\varepsilon^{\frac{1}{n+3}}.
\end{equation}
\end{theorem}

\textbf{Proof:} We prove it by contradiction. Suppose the conclusion fails, then there exist $\varepsilon_0$,
$s_0\in (-1,0)$ and $x_0\in \Sigma_{s_0}^{\infty}\cap B_{\frac{1}{2}}(0)$, such that
$|A|(x_0,s_0)>\left(\sqrt{\frac{1}{\kappa_0}}+C_2+C_3+C_2C_4\right)\varepsilon_0^{\frac{1}{n+3}},$ and
$\int^{0}_{-1}\int_{\Sigma^{\infty}_{s}\cap B_{1}(0)}|A|^2d\mu_s^{\infty}ds\leq \varepsilon_0$ .
Choosing $\delta=\varepsilon_0^{\frac{1}{n+3}}$, we know from (\ref{e3.3}) that
\begin{equation}\label{e3.7}
|A|(x,s_0) > |A|(x_0,s_0)-C_2 \delta>0 \ \ \ on \ \ B_{s_0}(x_0,\delta)\subset \Sigma^{\infty}_{s_0}\cap B_1(0).
\end{equation}
By (\ref{e3.4}), we know that there exists an interval $[a,a+\delta]\subset [-1,0]$ such that
$s_0 \in [a,a+\delta]$, and for any $s\in [a,a+\delta]$, $x\in \ B_{s_0}(x_0,\delta)$
\begin{equation}\label{e3.8}
|A|(x,s) \geq |A|(x,s_0)-C_3 \delta>0.
\end{equation}

By (\ref{e3.5}), we know that the metrics on the surfaces $\Sigma^{\infty}_{s}$ are all equivalent. In particular, we have
$B_{s}(x_0,\delta)\subset B_{s_0}(x_0,(1+C_4)\delta)$ for $s\in [a,a+\delta]$.
Thus we see from (\ref{e3.7}) and (\ref{e3.8}) that we have the estimate
\begin{equation}\label{e3.9}
|A|(x,s) > |A|(x_0,s_0)-(C_2+C_3+C_2C_4) \delta>0,
\end{equation}
for $s\in [a,a+\delta]$ and $x\in B_{s}(x_0,\delta)$. Therefore, we have by our assumption and Lemma \ref{lem3.2}
\begin{eqnarray*}
\varepsilon_0 & \geq & \int^{0}_{-1}\int_{\Sigma^{\infty}_{s}\cap B_1(0)}|A|^2d\mu_s^{\infty}ds
              \geq \int^{a+\delta}_{a}\int_{B_s(x_0,\delta)}|A|^{2}d\mu^{\infty}_sds\nonumber\\
            & \geq & (|A(x_{0},s_0)-(C_2+C_3+C_2C_4)\delta|)^{2} \int^{a+\delta}_{a}Vol_{s}(B_s(x_0,\delta))ds\nonumber\\
            & > & \frac{\delta^2}{\kappa_0}\delta\kappa_0 \delta^{n}=\delta^{n+3}=\varepsilon_0,
\end{eqnarray*}
which is a contradiction. This proves the theorem.
\hfill Q.E.D.

\vspace{.2in}

From the proof of the above theorem it is easy to get the elliptic case. That means,
\begin{corollary}\label{c1}
For any $0<\varepsilon\leq r_0^{n+2}$,  if
$\int_{\Sigma^{\infty}_{s}\cap B_1(0)}|A|^2d\mu_s^{\infty}\leq \varepsilon$, then
we have for each $s\in [-1,0]$
\begin{equation}\label{e3.10}
\max_{\Sigma^{\infty}_{s}\cap B_1(0)}|A| \leq
\left(\sqrt{\frac{1}{\kappa_0}}+C_2\right)\varepsilon^{\frac{1}{n+2}}.
\end{equation}
\end{corollary}

\vspace{.2in}

\textbf{Proof of The Main Theorem(Continued):}
We choose $\varepsilon_0$ is so small that
the right hand of (\ref{e3.6}) is smaller than 1 with $\varepsilon=2\varepsilon_0$.
But note that $|A_{\infty}|(0,0)=1$, thus, by Theorem \ref{thm3.3}, we must have
\begin{equation*}
\int^{0}_{-1}\int_{\Sigma^{\infty}_{s}\cap B_1(0)}|A_{\infty}|^2d\mu_s^{\infty}ds\geq 2\varepsilon_0.
\end{equation*}
In particular, for $k>>1$ sufficiently large,
\begin{equation*}
\int^{0}_{-1}\int_{\Sigma^{k}_{s}\cap B_1(0)}|A_k|^2d\mu_s^{k}ds> \varepsilon_0.
\end{equation*}
Rewriting the above integral back to the original surface, we see that we have
\begin{equation}\label{e3.10}
\lambda_k^n\int^{t_k}_{t_k-\lambda_k^{-2}}\int_{\Sigma_{t}\cap B_{\lambda_k^{-1}}(X_0)}|A|^2d\mu_tdt> \varepsilon_0.
\end{equation}
Set $t_k=T-\rho^2$ and $t_k-\lambda_k^{-2}=T-(\rho')^2$, then $(\rho')^2-\rho^2=\lambda_k^{-2}$.
As $t_k\to T$ as $k\to \infty$, we see that as $k$ sufficiently large and fixed, we have
$0\leq \rho <\rho' \leq \frac{r}{2}$, and
\begin{equation*}
((\rho')^2-\rho^2)^{-\frac{n}{2}}\int^{T-\rho^2}_{T-(\rho')^2}\int_{\Sigma_t \cap B_{((\rho')^2-\rho^2)^{\frac{1}{2}}}}|A|^2d\mu_t dt> \varepsilon_0,
\end{equation*}
which contradicts our assumption. This proves the Main Theorem.
\hfill Q.E.D.

\vspace{.1in}

The equation (\ref{e3.10}) tells us the set of singularities along the mean curvature flow is defined by
$${\mathcal F}=\{X\in M|\exists \lambda_k\to \infty, t_k\to T \lim_{k\to \infty}
\lambda_k^{n}\int ^{t_k}_{t_k-\lambda_k^{-2}}\int_{\Sigma_t\cap
B_{\lambda_k^{-1}}(X)}|A|^2\geq\varepsilon_0\}
$$

If $\Sigma_0$ is a surface, i.e., $n=2$, as a corollary of the above theorem, we prove,

\begin{corollary}
${\cal H}^2({\mathcal F})=0$.
\end{corollary}
{\it Proof.} Fix $\delta>0$. From the definition of $\mathcal F$ we know that for
any $X\in\mathcal F$, there exists $r>0$ and $t<T$ such that $t-r^2>T-\delta$,
$r\leq\delta/10$ and $$
\int ^{t}_{t-r^{2}}\int_{\Sigma_t\cap
B_{r}(X)}|A|^2\geq\frac{\varepsilon_0}{2}r^{2}.
$$
 By Theorem 3.3 in \cite{Sim}, we can choose a disjoint
family of balls $\{B_{r_j}(X_j)\}$ with $X_j\in\mathcal F$ and
$r_j<\delta/10$ such that the family $\{B_{5r_j}(X_j)\}$ covers
$\mathcal F$ and for all $j\in \mathbf N$
$$\int_{t_j-r_j^2}^{t_j}\int_{\Sigma_t\cap
B_{r_j}(X_j)}|A|^2\geq r_j^2\varepsilon_0/2.
$$
We then estimate
\begin{eqnarray*}
{\cal H}^2_\delta({\mathcal F})&\leq& c\sum_{j=1}^\infty
r_j^2<c\varepsilon_0^{-1}\sum_{j=1}^\infty
\int_{t_j-r_j^2}^{t_j}\int_{\Sigma_t\cap B_{r_j}(X_j)}|A|^2\leq c\varepsilon_0^{-1}\sum_{j=1}^\infty
\int_{T-\delta}^{T}\int_{\Sigma_t\cap B_{r_j}(X_j)}|A|^2\\&\leq&
c\varepsilon_0^{-1}\int_{T-\delta}^T\sum_{j=1}^\infty\int_{\Sigma_t\cap
B_{r_j}(X_j)}|A|^2\leq
c\varepsilon_0^{-1}\int_{T-\delta}^T\int_{\Sigma_t}|A|^2.
\end{eqnarray*}
By the Gauss equation
$$
R_{1212}=K_{1212}+(h_{11}^\alpha h_{22}^\alpha - h_{12}^\alpha
h_{12}^\alpha),
$$
we get
$$
|A|^2=|H|^2-2R_{1212}+2K_{1212},
$$
were $K$ is the curvature of $M$ and $R$ is the curvature of $\Sigma$. We therefore have by the Gauss-Bonnet formula that,
$$
\int_{\Sigma_t}|A|^2d\mu_t\leq \int_{\Sigma_t}|H|^2d\mu_t
+C\mu_t(\Sigma_t)+8\pi (g-1),
$$
where $g$ is the genus of the initial surface $\Sigma_0$. Since
$$
\frac{\partial}{\partial
t}\int_{\Sigma_t}d\mu_t=-\int_{\Sigma_t}|H|^2 d\mu_t,
$$ we have
\begin{equation}\label{eq18}\mu_t(\Sigma_t)\leq \mu_0(\Sigma_0)~{\rm
and}~\int_0^T\int_{\Sigma_t}|H|^2d\mu_tdt \leq \mu_0(\Sigma_0)
.\end{equation} So,
\begin{equation}\label{eq19}
\int_{\Sigma_t}|A|^2d\mu_t\leq \int_{\Sigma_t}|H|^2d\mu_t +C,
\end{equation} and consequently,
$$\int^T_0 \int_{\Sigma_t}|A|^2 d\mu_t dt\leq C,
$$
and therefore,
$$
\lim_{\delta\to 0}\int_{T-\delta}^T\int_{\Sigma_t}|A|^2=0. $$
Thus, ${\cal H}^2 ({\mathcal F})=\lim_{\delta\to 0}{\cal
H}^2_\delta({\mathcal F})=0$. We finish the proof of our corollary.
 \hfill Q. E. D.


{\small\begin{thebibliography}{999}

\bibitem {CY} B. Chen and L. Yin, {\em Uniqueness and pseudolocality theorems of the mean curvature flow},
Comm. Anal. Geom.,  {\bf 15} (2007), no.3, 435-490.

\bibitem {CL1} J. Chen and J. Li, {\em Mean curvature flow of
surface in 4-manifolds}, Adv. Math., {\bf 163} (2001), 287-309.

\bibitem {CLT} J. Chen, J. Li and G. Tian {\em Two-dimensional graphs moving by mean curvature flow},
Acta Math. Sin., {\bf 18} (2002), 209-224.

\bibitem {Ec}
K. Ecker, {\em On regularity for mean curvature flow of hypersurfaces}, Calc. Var., {\bf 3} (1995), 107-126.

\bibitem {EH} K. Ecker and G. Huisken, {\em Mean curvature evolution of
entire graphs}, Ann. of Math., {\bf 130} (1989), 453-471.

\bibitem {ES}
J. Eells and J. H. Sampson, {\em Harmonic maps of Riemannian manifolds}, Amer. J. Math., {\bf 86} (1964), 109-160.

\bibitem {HanS2}
X. Han and J. Sun, {\em $\varepsilon_{0}$-regularity for mean
curvature flow from surface to flat Riemannian manifold}, to appear
in Acta Math. Sin..

\bibitem{Hu1}
G. Huisken, {\em Local and global behaviour of hypersurfaces moving
by mean curvature}, Proceedings of Symposia in Pure Mathematics,
{\bf 54} (1993), Part I, 175-191 .

\bibitem {Il}
T. Iimanen, {\em Singularities of mean curvature flow of surfaces}, preprint.

\bibitem {SchU} R. Schoen and K. Uhlenbeck, {\em A regularity theory of harmonic maps}, J.
Diff. Geom., {\bf 17} (1982), 307-335 and {\bf 18} (1983), 329.

\bibitem {Sim} L. Simon, {\em Lectures on Geometric Measure Theory}, Proc. of the CMA, Vol. 3 (1983).

\bibitem {Wa1} M.-T. Wang, {\em Mean curvature flow of surfaces in
Einstein four manifolds}, J. Diff. Geom., {\bf 57} (2001), 301-338.

\bibitem {Wa2} M.-T. Wang, {\em Long time existence and convergence of graphic mean curvature flow in
arbitray codimension}, Invent. Math., {\bf 148} (2002), 525-543.

\end{thebibliography}}
\end{document}